\documentclass{article}
\usepackage[sort, numbers]{natbib}
\usepackage{amsmath,amsfonts,amsthm} % Math packages
\usepackage{times, fullpage,color}
\usepackage{mathtools,url,tikz}
\newcommand{\norm}[1]{{\left\lVert#1\right\rVert}}

\newtheorem{theorem}{Theorem}[section]
\newtheorem{definition}[theorem]{Definition}
\newtheorem{corollary}[theorem]{Corollary}
\newtheorem{example}[theorem]{Example}
\newtheorem{lemma}[theorem]{Lemma}

\newcommand{\modAT}[1]{{\color{black}#1}}
\newcommand{\modFG}[1]{{\color{black}#1}}
\newcommand{\modFGR}[1]{{\color{black}#1}}
\newcommand{\modATR}[1]{{\color{black}#1}}

\newcommand{\inner}[2]{{#1^{\tsp} #2}}
\def\tsp{{\sf T}}

\title{On the  worst-case complexity of the gradient method \\with exact line search for smooth strongly convex functions}
\author{Etienne de Klerk\thanks{Tilburg University and Delft University of Technology, The Netherlands,  \texttt{E.deKlerk@uvt.nl}} \and Fran\c{c}ois Glineur\thanks{
UCL / CORE and ICTEAM, Louvain-la-Neuve, Belgium, \texttt{Francois.Glineur@uclouvain.be}, \texttt{Adrien.Taylor@uclouvain.be}. A.B. Taylor is a F.R.I.A. fellow (F.R.S.-FNRS).
The UCL/CORE authors are supported by the Belgian Interuniversity Attraction Poles, and by the
ARC grant 13/18-054 (Communaut\'e fran\c{c}aise de Belgique).
}
 \and  Adrien B. Taylor\footnotemark[2]}

\begin{document}

\maketitle

\begin{abstract}
We consider the gradient (or steepest) descent method with exact line {search applied} to a strongly convex function
with Lipschitz continuous gradient. We establish the exact worst-case rate of convergence {of this scheme}, and show that this
worst-case behavior is exhibited by a certain convex quadratic function. We also give \modFG{the} \modATR{tight} worst-case complexity bound for a noisy variant of gradient descent method,
where exact line-search is performed in a search direction that differs from negative gradient by at most a prescribed relative tolerance.

The proof\modATR{s are} computer-assisted, and rel\modATR{y} on the resolution\modATR{s} of semidefinite programming performance estimation problems as introduced in the paper [Y. Drori and M. Teboulle.
 Performance of first-order methods for smooth convex minimization: a novel approach. \emph{Mathematical Programming}, 145(1-2):451-482, 2014].
\end{abstract}


\vspace{0.2cm}


\section{Introduction}

The gradient (or steepest) descent method for unconstrained method was devised by Augustin-Louis Cauchy (1789--1857) in the 19th century, and remains
one of the most iconic algorithms for unconstrained optimization. Indeed, it is usually the first algorithm that is taught during {introductory courses} on nonlinear optimization.
It is therefore somewhat surprising that the worst-case converge{nce} rate of the method is not yet precisely understood for smooth strongly convex functions.

In this paper, we settle the worst-case convergence rate question of the gradient descent method with exact {line search}
for strongly convex, continuously differentiable functions $f$ with Lipschitz continuous gradient.
Formally we consider the following function class.

\begin{definition}
A continuously differentiable function $f:\mathbb{R}^n \rightarrow \mathbb{R}$ is called $L$-smooth, $\mu$-strongly
convex with parameters $L > 0$ and $\mu >0$ if
\begin{enumerate}
\item
$\mathbf{x} \mapsto f(\mathbf{x}) -  \frac{\mu}{2}\|\mathbf{x}\|^2$ is a convex function on $\mathbb{R}^n$, where the norm is the Euclidean norm;
\item
$\|\nabla f(\mathbf{x}+\Delta \mathbf{x}) - \nabla f(\mathbf{x}) \| \le L \|\Delta \mathbf{x}\|$ holds for all $\mathbf{x} \in \mathbb{R}^n$ and $\Delta \mathbf{x} \in \mathbb{R}^n$.
\end{enumerate}
The class of $L$-smooth, $\mu$-strongly convex functions \modFG{on $\mathbb{R}^n$} will be denoted by
$\mathcal{F}_{\mu,L}(\mathbb{R}^n)$.
\end{definition}
Note that, if  $f$ is twice continuously differentiable, then $f \in \mathcal{F}_{\mu,L}(\mathbb{R}^n)$ is
equivalent to
\[ L I \succeq \nabla^2 f(\mathbf{x}) \succeq \mu I \quad \forall \mathbf{x} \in \mathbb{R}^n
\]
where the notation $A \succeq B$ for symmetric matrices $A$ and $B$ means the matrix $A-B$ is positive semidefinite, and $I$ is the identity matrix.
Equivalently,  the eigenvalues of the Hessian matrix $\nabla^2 f(\mathbf{x})$ lie in the interval $[\mu,L]$ for all $\mathbf{x}$.

The gradient method with exact line search may be described as follows.

\begin{center}
{\small \fbox{
\parbox{0.7\textwidth}{
        \textbf{Gradient descent method with exact line search}
  \begin{itemize}
  \item[] \textbf{Input:} $f\in\mathcal{F}_{\mu,L}(\mathbb{R}^n)$, $\mathbf{x}_0\in \mathbb{R}^n$.\\[-0.2cm]
  \item[] \textbf{for} $i=0,1,\ldots$ \\ [-0.2cm]
  \item[]\hspace{1cm} {$\gamma=\text{argmin}_{\gamma \in \mathbb{R}}f\left(\mathbf{x}_i-\gamma \nabla f(\mathbf{x}_i)\right)$} \\ [-0.2cm]
  \item[]
   \hspace{1cm} {$\mathbf{x}_{i+1}=\mathbf{x}_i-\gamma \nabla f(\mathbf{x}_i)$} \\ [-0.2cm]
  \end{itemize}
       }}}
\end{center}

Our main result may now be stated concisely.

\begin{theorem}
\label{th:main}
Let $f\in \mathcal{F}_{\mu,L}(\mathbb{R}^n)$, $\mathbf{x}_*$ a global minimizer of $f$ on $\mathbb{R}^n$, and $f_* = f(\mathbf{x}_*)$.
Each iteration of the gradient method with exact line search satisfies
\begin{equation}
\label{eq:bound}
{f(\mathbf{x}_{i+1})} - f_*\leq \left(\frac{L-\mu}{L+\mu}\right)^2 \left({f(\mathbf{x}_i)}-f_*\right) \quad i = 0,1,\ldots
\end{equation}
\end{theorem}
Note that the result in Theorem \ref{th:main}\modFG{, which establises a global linear convergence rate on objective function accuracy,} is known for the case of quadratic functions in $\mathcal{F}_{\mu,L}(\mathbb{R}^n)$, that
is for functions of the form
\[
f(\mathbf{x}) = \frac{1}{2}\mathbf{x}^{\tsp}Q\mathbf{x} + \mathbf{c}^{\tsp}\mathbf{x}
\]
where $\mathbf{c} \in \mathbb{R}^n$, and the eigenvalues of the $n\times n$ symmetric positive definite matrix $Q$ lie in the interval $[\mu,L]$;
see e.g.\ \cite[\S1.3]{bertsekas1999nonlinear}, \cite[pp. 60--62]{polyak1987intro}, or \cite[pp. 235--238]{luenberger2008linear}.
Moreover, the bound \eqref{eq:bound} is known to be tight  for the
following example.

\begin{figure}
\begin{center}
\begin{tikzpicture}[yscale=1,xscale=1]
\tikzset{fleche/.style={->, >=latex, thick},
	    flecheB/.style={<->, >=latex, very thick},
 	    droite/.style={thick, dashed}};
 	
\pgfmathsetmacro{\muExB}{.2};
\pgfmathsetmacro{\LExB}{2};
\pgfmathsetmacro{\rhoExB}{(\LExB-\muExB)/(\LExB+\muExB)};
\pgfmathsetmacro{\fxs}{(1/\LExB+1/\muExB)/2};
\pgfmathsetmacro{\lAxis}{sqrt(2*\fxs/\muExB)}
\pgfmathsetmacro{\sAxis}{sqrt(2*\fxs/\LExB)};

\pgfmathsetmacro{\nbEllipses}{8};%Number of ellipses to show

\foreach \i in {1,...,\nbEllipses}
{
\draw [dashed, color=black, domain=-0:360, samples=50] plot({cos(\x)*\lAxis*\rhoExB^(\i-1)}, {sin(\x)*\sAxis*\rhoExB^(\i-1)});
}

\node (XS) at (0,0) {$\bullet$};
\draw (XS) node[above left] {$\mathbf{x}_*$};
\node (X0) at (1/\muExB,1/\LExB) {$\bullet$};
\draw (X0) node[right ] {${\mathbf{x}_0} = [1/\mu,1/L]^\tsp$};
\node (X1) at (\rhoExB/\muExB,-\rhoExB/\LExB) {$\bullet$};
\draw (X1) node[below ] {$\mathbf{x}_1$};
\node (X2) at (\rhoExB^2/\muExB,\rhoExB^2/\LExB) {$\bullet$};
\draw (X2) node[above ] {$\mathbf{x}_2$};
\node (X3) at (\rhoExB^3/\muExB,-\rhoExB^3/\LExB) {$\bullet$};
\draw (X3) node[below ] {$\mathbf{x}_3$};
\node (X4) at (\rhoExB^4/\muExB,\rhoExB^4/\LExB) {$\bullet$};
\draw (X4) node[above ] {$\mathbf{x}_4$};
\node (X5) at (\rhoExB^5/\muExB,-\rhoExB^5/\LExB) {$\bullet$};
\draw (X5) node[below ] {$\mathbf{x}_5$};
\node (X6) at (\rhoExB^6/\muExB,\rhoExB^6/\LExB) {$\bullet$};
\draw (X6) node[above ] {$\mathbf{x}_6$};
\node (X7) at (\rhoExB^7/\muExB,-\rhoExB^7/\LExB) {$\bullet$};
\draw (X7) node[below ] {$\mathbf{x}_7$};
\draw[thick] (X0.center)--(X1.center)--(X2.center)--(X3.center)--(X4.center)--(X5.center) -- (X6.center) -- (X7.center);
\draw[fleche] (X0.center)--(1/\muExB-.2*1/\muExB/3,1/\LExB-2*1/\LExB/3);
\draw[fleche] (X1.center)--(\rhoExB*1/\muExB-.2*1/\muExB/3,-\rhoExB*1/\LExB+2*1/\LExB/3);
\draw [fleche] (-\lAxis-.5,0) -- (\lAxis+1,0);
\draw [fleche] (0,-\sAxis-.5) -- (0,\sAxis+1);
\draw (\lAxis+1,0) node[ right] {$x_1$};
\draw (0,\sAxis+1) node[above right] {$x_2$};
\draw (0,\sAxis) node[above right] {$\frac{1}{\sqrt{L}}$};
\draw (\lAxis,0) node[below right] {$\frac{1}{\sqrt{\mu}}$};
\end{tikzpicture}
\caption{Illustration of Example \ref{eg:quad example} for the case $n=2$ (small arrows indicate direction of negative gradient). \label{fig:quad_ex}}
\end{center}
\end{figure}
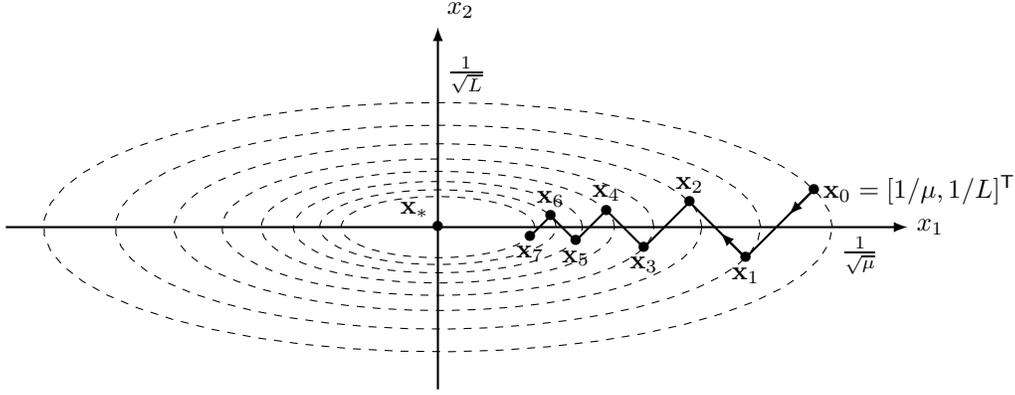

\begin{example}
\label{eg:quad example}
Consider the following quadratic function from \cite[Example on p. 69]{bertsekas1999nonlinear}:
\[
f(\mathbf{x}) = \frac{1}{2}\sum_{i=1}^n \lambda_i x_i^2
\]
where
\[
0 < \mu = \lambda_1 \le \lambda_2 \le \ldots \le \lambda_{n} = L,
\]
 and the starting point
\[
\mathbf{x}_0 = (\frac{1}{\mu}, 0, \ldots, 0, \frac{1}{L})^\tsp.
\]
One may readily \modFG{check that the gradient at $\mathbf{x}_0$ is equal to
\[
\nabla f(\mathbf{x}_0) = (1, 0, \ldots, 0, 1)^\tsp
\]
and that the minimum of the line-search from $\mathbf{x}_0$ in that direction
is attained for step $\gamma = \frac{2}{L+\mu}$. One therefore obtains}
\[
\mathbf{x}_1 = \left(\frac{L-\mu}{L+\mu}\right)(1/\mu, 0, \ldots, 0, -1/L)^{\tsp},
\]
and, for all $i=0,1,\ldots$
\[
\mathbf{x}_{2i} = \left(\frac{L-\mu}{L+\mu}\right)^{2i}x_0, \;\;\; \mathbf{x}_{2i+1} = \left(\frac{L-\mu}{L+\mu}\right)^{2i}x_1.
\]
Since $f_* = 0$, it is straightforward to verify that equality
\[
{f(\mathbf{x}_{i+1})} - f_* = \left(\frac{L-\mu}{L+\mu}\right)^2 ({f(\mathbf{x}_i)}-f_*) \quad i = 0,1,\ldots,
\]
holds as required. \qed
\end{example}
The construction in Example \ref{eg:quad example} is illustrated in Figure \ref{fig:quad_ex} \modFG{in the case $n=2$, where the  ellipses shown are level curves of the objective function. Each step from $\mathbf{x}_i$ to $\mathbf{x}_{i+1}$  is orthogonal to the ellipse at $\mathbf{x}_i$  (since it uses the steepest descent direction) and tangent to the ellipse at $\mathbf{x}_{i+1}$ (because of the exact line-search direction), hence successive steps are orthogonal to each other.}

As an immediate consequence of Theorem \ref{th:main} and  Example  \ref{eg:quad example}, one has the following tight bound on the number of steps needed to obtain $\epsilon$-relative accuracy on the objective function for a given $\epsilon > 0$.

\begin{corollary}
Given $\epsilon > 0$, the gradient method with exact line search yields
a solution with relative accuracy $\epsilon$ for any function $f\in \mathcal{F}_{\mu,L}(\mathbb{R}^n)$ after at most $N = \left\lceil \frac{1}{2}\log \left( \frac{1}{\epsilon}\right) / \log\left(\frac{L+\mu}{L-\mu}\right)\right\rceil$ iterations, i.e.
\[
\frac{f(\mathbf{x}_N) - f_*}{f(\mathbf{x}_0) - f_*} \le \epsilon,
\]
where $\mathbf{x}_0$ is the starting point.
Moreover, this iteration bound is tight for the quadratic function defined in {Example~\ref{eg:quad example}}.
\end{corollary}

For non-quadratic functions in $\mathcal{F}_{\mu,L}(\mathbb{R}^n)$, only bounds weaker than \eqref{eq:bound} are known.
For example, in \cite[p. 240]{luenberger2008linear}, the following bound is shown:
\[
{(f(\mathbf{x}_{i+1})-f_*)} \le \left(1 - \frac{\mu}{L} \right){(f(\mathbf{x}_{i})-f_*)} \quad i = 0,1,\ldots
\]

In  \cite[Theorem 3.4]{nocedal2006numerical} a stronger result than Theorem \ref{th:main} was claimed, but this was retracted
in a subsequent erratum\footnote{The erratum is available at: \url{http://users.iems.northwestern.edu/~nocedal/book/2ndprint.pdf}}, and only an asymptotic result is claimed in the erratum.

A result related to Theorem \ref{th:main} is given in \cite{nemirovskinotes1999} where
Armijo-rule line search is used instead of exact line search. An explicit rate in the strongly convex case is given there
 in Proposition 3.3.5 on page 53 (definition of the method is (3.1.2) on page 44).
 More general upper bounds on the convergence rates of gradient-type methods for convex functions may be found in
  the books \cite{nemirovsky1983problem,nesterov2004intro}.
We mention one more particular result by Nesterov \cite{nesterov2004intro} that is similar to our main result in Theorem \ref{th:main}, \modFG{but that uses a fixed step-length and relies on the initial distance to the solution}.
\begin{theorem}[Theorem 2.1.15 in \cite{nesterov2004intro}] \label{thNesterov}
Given $f\in \mathcal{F}_{\mu,L}(\mathbb{R}^n)$ and $\mathbf{x}_0 \in \mathbb{R}^n$, the gradient descent method with {fixed step length} $\modFG{\gamma}=\frac{2}{\mu+L}$ generate iterates
$\mathbf{x}_i$ $(i=0,\modFG{1,}2,\ldots)$ that satisfy
\[
{f(\mathbf{x}_{i})} - f_* \leq \frac{L}{2}\left(\frac{L-\mu}{L+\mu}\right)^{2i} \left\|\mathbf{x}_0 - \mathbf{x}_* \right\|^2 \quad i = 0,1,\ldots
\]
\end{theorem}
Note that this result does not imply Theorem \ref{th:main}.

\section{Background results}
In this section we collect some known results on strongly convex functions and on the gradient method.
We will need these results in the proof of our main result, Theorem \ref{th:main}.
%To simplify our notation we will denote
%\[
%f_i = f(x_i), \; g_i = \nabla f(x_i), \; i=1,2,\ldots
%\]

\subsection{Properties of the gradient method with exact line search}
{Let} $\mathbf{x}_i$ ($i=1,2,\ldots,N$) {be the} iterates produced by the gradient method with exact line search {started at $\mathbf{x}_0$}. Those iterates are defined by the following two conditions for $i =0,1,\ldots,N-1$
\begin{align}
\label{eq:negative gradient} \mathbf{x}_{i+1}-\mathbf{x}_i+\gamma \nabla f(\mathbf{x}_i)&=0, \ \text{for some } \gamma\geq 0,\\ \label{eq:exact line search}
\nabla f(\mathbf{x}_{i+1})^{\tsp} (\mathbf{x}_{i+1} - \mathbf{x}_i) &= 0
\end{align}
where the first condition \eqref{eq:negative gradient} states that we move in the direction of the negative gradient, and the second condition \eqref{eq:exact line search} expresses the exact line search condition.

A consequence of those conditions is that successive gradients are orthogonal, i.e.
\begin{equation}
\label{eq:successive gradients orthogonal}
\nabla f(\mathbf{x}_{i+1})^{\tsp} \nabla f(\mathbf{x}_{i}) = 0 \quad i = 0,1,\ldots,N-1.
\end{equation}

Instead of relying on conditions \eqref{eq:negative gradient}--\eqref{eq:exact line search} that define the iterates of the gradient method with exact line search, our analysis will be based on the weaker conditions \eqref{eq:exact line search}--\eqref{eq:successive gradients orthogonal}, which are also satisfied by other sequences of iterates.

\subsection{Interpolation with functions in $\mathcal{F}_{\mu,L}(\mathbb{R}^n)$}
We \modFG{now} consider the following interpolation problem over the class of functions $\mathcal{F}_{\mu,L}(\mathbb{R}^n)$.
\begin{definition}Consider an integer $N \ge 1$ and given data
$\{{(\mathbf{x}_i,f_i,\mathbf{g}_i)}\}_{i \in \{0,1,\ldots,N\}}$ where
$\mathbf{x}_i \in \mathbb{R}^n$, $f_i\in \mathbb{R}$ and $\mathbf{g}_i\in \mathbb{R}^n$.
If there exists a function $f \in  \mathcal{F}_{\mu,L}(\mathbb{R}^n)$
such that
\[
f(\mathbf{x}_i) = f_i, \; \nabla f(\mathbf{x}_i) = \mathbf{g}_i, \;\forall i \in \{0,1,\ldots,N\},
\]
then we say that $\{{(\mathbf{x}_i,f_i,\mathbf{g}_i)}\}_{i \in \{0,1,\ldots,N\}}$ is $\mathcal{F}_{\mu,L}$-interpolable.
\end{definition}

A necessary and sufficient condition for $\mathcal{F}_{\mu,L}$-interpolability in given in the next theorem, taken from \cite{PEP2016smooth}.

\begin{theorem}[\cite{PEP2016smooth}]
\label{th:interpolation}
A data set  $\{{(\mathbf{x}_i,f_i,\mathbf{g}_i)}\}_{i \in \{0,1,\ldots,N\}}$ is $\mathcal{F}_{\mu,L}$-interpolable if and only if the following inequality
\[
f_i - f_j - \inner{\mathbf{g}_j}{(\mathbf{x}_i-\mathbf{x}_j)} \geq \frac{1}{2(1-\mu/L)}\left( \frac{1}{L}\norm{\mathbf{g}_i-\mathbf{g}_j}^2+ \mu \norm{\mathbf{x}_i-\mathbf{x}_j}^2 - 2\frac{\mu}{L} \inner{(\mathbf{g}_j-\mathbf{g}_i)}{(\mathbf{x}_j-\mathbf{x}_i)}\right)
\]
holds for all ${i\neq j} \in \{0,1,\ldots,N\}$.
\end{theorem}
In principle, Theorem \ref{th:interpolation} allows one to generate all possible valid inequalities that hold for
 functions in $\mathcal{F}_{\mu,L}(\mathbb{R}^n)$
in terms of their function values and gradients at a  set of points $\mathbf{x}_0,\ldots,\mathbf{x}_N$. This will be the essential for the proof of our main result, Theorem \ref{th:main}.

\section{A performance estimation problem}
The proof technique we will use for Theorem \ref{th:main}
is inspired by recent work on the so-called performance estimation problem,
as introduced in \cite{drori2014} and further developed in \cite{PEP2016smooth}.
The idea is to formulate the computation of the worst-case behavior of certain
iterative methods as an explicit semidefinite programming (SDP) problem.
We first recall the definition of SDP problems (in a form that is suitable to our purposes).

\subsection{Semidefinite programs}
We will consider semidefinite programs (SDPs) of the form
\begin{equation}
\label{eq:primal sdp}
\max_{X = (x_{ij}) \in \mathbb{S}^n, X \succeq 0, \mathbf{u} \in \mathbb{R}^{\ell}} \left\{ \sum_{i,j = 1}^n c_{ij}x_{ij} + \mathbf{c}^\tsp \mathbf{u}\; \left| \; \sum_{i,j = 1}^n a^{(k)}_{ij}x_{ij} + \mathbf{a}_k^\tsp \mathbf{u} \le b_k \quad k = 1,\ldots,m \right\},\right.
\end{equation}
where $\mathbb{S}^n$ is the set of symmetric matrices of size $n$, and matrices $A_k = \left(a^{(k)}_{ij}\right) \in \mathbb{S}^n$ and the matrix $C = (c_{ij}) \in \mathbb{S}^n$ are given, as well
as the scalars $b_k$ and vectors $\mathbf{a}_k \in \mathbb{R}^\ell$ ($k = 1,\ldots,m$), and $\mathbf{c} \in \mathbb{R}^\ell$.

Since every positive semidefinite matrix $X \in \mathbb{S}^n$ is a Gram matrix, there exist vectors $\mathbf{v}_1,\ldots,\mathbf{v}_n \in \mathbb{R}^n$
such that $x_{ij} = \mathbf{v}_i^{\tsp}\mathbf{v}_j$ for all $i,j$.
Thus the SDP problem \eqref{eq:primal sdp} may be equivalently rewritten as
\begin{equation}
\label{eq:primal sdp Gram}
\max_{\mathbf{v}_i \in \mathbb{R}^n, \mathbf{u} \in \mathbb{R}^{\ell}} \left\{ \sum_{i,j = 1}^n c_{ij}\mathbf{v}_i^{\tsp}\mathbf{v}_j + \mathbf{c}^\tsp \mathbf{u}\; \left| \; \sum_{i,j = 1}^n a^{(k)}_{ij}\mathbf{v}_i^{\tsp}\mathbf{v}_j + \mathbf{a}_k^\tsp \mathbf{u} \le b_k \quad k = 1,\ldots,m\right\}\right.
\end{equation}
which features terms that are linear in the inner products $\mathbf{v}_i^{\tsp}\mathbf{v}_j$ in the objective function and constraints.
The associated dual SDP problem is
\begin{equation}
\label{eq:dual SDP}
\min_{\mathbf{y} \in \mathbb{R}^m, \mathbf{y} \ge \mathbf{0}} \left\{\mathbf{b}^{\tsp}\mathbf{y} \; \left| \; \sum_{k=1}^m y_kA_k - C \succeq 0, \; \sum_{k=1}^m y_k\mathbf{a}_k   = \mathbf{c}\right\}.\right.
\end{equation}
We will later use the fact that each dual variable $y_k$ may be viewed as a (Lagrange) multiplier of the primal constraint {$\sum_{i,j = 1}^n a^{(k)}_{ij}\mathbf{v}_i^{\tsp}\mathbf{v}_j + \mathbf{a}_k^\tsp \mathbf{u}\le b_k$}.

%By the weak duality theorem, for all feasible primal-dual solutions $X$ and $y$ (resp.) one has
%\begin{equation}
%\label{eq:weak duality}
%0 \le \sum_{i,j = 1}^n c_{ij}x_{ij} - \sum_{k=1}^m \sum_{i,j = 1}^n y_k a^{(k)}_{ij}x_{ij}.
%\end{equation}

\subsection{Performance estimation of the gradient method with exact line search}
Consider the following SDP problem, for fixed parameters $N \ge 1$, $R>0$, $\mu> 0$ and $L>\mu$:

\begin{equation}
\label{eq:pep sdp}
\left.
\begin{array}{llcl}
& \max f_N - f_* & & \\
\mbox{subject to} &&& \\
&\mathbf{g}_{i+1}^{\tsp} (\mathbf{x}_{i+1} - \mathbf{x}_i) &=& 0 \quad i \in {\{0,1,\ldots,N-1\} } \\
&\mathbf{g}_{i+1}^{\tsp} \mathbf{g}_i &= &0 \quad i \in {\{0,1,\ldots,N-1\} } \\
&\{{(\mathbf{x}_i,f_i,\mathbf{g}_i)}\}_{i \in {\{*,0,1,\ldots,N\}}} & \text{is} & \mbox{ $\mathcal{F}_{\mu,L}$-interpolable} \\
&\mathbf{g}_* &=& 0 \\
%&f_* &=& 0 \\
&f_0-f_*&\le&R,\\ \\
\end{array}
\right\}
\end{equation}
where the variables are $\mathbf{x}_i \in \mathbb{R}^n$, $f_i\in \mathbb{R}$ and $\mathbf{g}_i\in \mathbb{R}^n$ ($i \in \{*,0,1,\ldots,N\}$).

Note that this is indeed an SDP problem of the form \eqref{eq:primal sdp Gram}, with dual problem of the form \eqref{eq:dual SDP}, since equalities and interpolability conditions are linear in the inner products of variables $\mathbf{x}_i$ and $\mathbf{g}_i$.
\begin{lemma}
\label{lemma:pep}
The optimal value of the above SDP problem \eqref{eq:pep sdp} is an upper bound on $f(\mathbf{x}_N) - f_*$, where $f$ is any function from $\mathcal{F}_{\mu,L}(\mathbb{R}^n)$, $f_*$ is its minimum and $\mathbf{x}_N$ is the $N$th iterate of the gradient method with
exact line search applied to $f$ from any starting point $\mathbf{x}_0$ that satisfies $f(\mathbf{x}_0) {-f_*}\le R$.
\end{lemma}
\proof
Fix any $f \in \mathcal{F}_{\mu,L}(\mathbb{R}^n)$, and let $\mathbf{x}_0,\ldots,\mathbf{x}_N$ be the iterates of the gradient method with exact line search applied to $f$.
Now a feasible solution to the SDP problem is given by
\[
\mathbf{x}_i, \; f_i = f(\mathbf{x}_i), \; \mathbf{g}_i = \nabla f(\mathbf{x}_i) \;\quad i \in \{*,0,\ldots,N\}.
\]
The objective function value at this feasible point is $f_N = f(\mathbf{x}_N)$, so that the optimal value of the SDP is an upper bound on
$f(\mathbf{x}_N)-f_*$. %(Note the assumption $f_* = 0$ is without loss of generality.
 \qed

We are now ready to give a proof of our main result. We already mention that the SDP relaxation \eqref{eq:pep sdp} is not used directly in the proof, but was used
to devise the proof, in a sense that will be explained later.

\section{Proof of Theorem \ref{th:main}}
A little reflection shows that, to prove Theorem \ref{th:main}, we need only consider one iteration of the gradient method with exact {line search}.
Thus we consider only the first iterate, given by $\mathbf{x}_0$ and $\mathbf{x}_1$, as well as the minimizer $\mathbf{x}_*$ of $f \in \mathcal{F}_{\mu,L}$.

Set $f_i = f(\mathbf{x}_i)$ and $\mathbf{g}_i = \nabla f(\mathbf{x}_i)$ for $i \in \{*,0,1\}$. Note that $\mathbf{g}_* = \mathbf{0}$.
The following five inequalities are now satisfied:
\begin{eqnarray*}
1: &&f_{0}\geq f_1 + \inner{\mathbf{g}_1}{(\mathbf{x}_{0}-\mathbf{x}_1)} +\frac{1}{2(1-\mu/L)}\left( \frac{1}{L}\norm{\mathbf{g}_0-\mathbf{g}_1}^2+ \mu \norm{\mathbf{x}_0-\mathbf{x}_1}^2 - 2\frac{\mu}{L} \inner{(\mathbf{g}_1-\mathbf{g}_0)}{(\mathbf{x}_1-\mathbf{x}_0)}\right)\\
2: &&f_*\geq f_0 + \inner{\mathbf{g}_0}{(\mathbf{x}_*-\mathbf{x}_0)} +\frac{1}{2(1-\mu/L)}\left( \frac{1}{L}\norm{\mathbf{g}_*-\mathbf{g}_0}^2+ \mu \norm{\mathbf{x}_*-\mathbf{x}_0}^2 - 2\frac{\mu}{L} \inner{(\mathbf{g}_0-\mathbf{g}_*)}{(\mathbf{x}_0-\mathbf{x}_*)}\right)\\
3: &&f_*\geq f_1 + \inner{\mathbf{g}_1}{(\mathbf{x}_*-\mathbf{x}_1)} +\frac{1}{2(1-\mu/L)}\left( \frac{1}{L}\norm{\mathbf{g}_*-\mathbf{g}_1}^2+ \mu \norm{\mathbf{x}_*-\mathbf{x}_1}^2 - 2\frac{\mu}{L} \inner{(\mathbf{g}_1-\mathbf{g}_*)}{(\mathbf{x}_1-\mathbf{x}_*)}\right)\\
4: &&-\mathbf{g}_{0}^{\tsp} \mathbf{g}_{1}\geq 0\\
5: &&\mathbf{g}_{1}^{\tsp} (\mathbf{x}_0-\mathbf{x}_1)\geq 0.\\
\end{eqnarray*}
Indeed, the first three inequalities are the $\mathcal{F}_{\mu,L}$-interpolability conditions,
the fourth inequality is a relaxation of \eqref{eq:successive gradients orthogonal}, and the fifth inequality is a relaxation
of \eqref{eq:exact line search}.

We aggregate these five inequalities by
defining the following positive multipliers,
\begin{equation}
\label{eq:multipliers}
 y_1=\frac{L-\mu}{L+\mu}, \quad y_2=2\mu\frac{(L-\mu)}{(L+\mu)^2}, \quad y_3=\frac{2\mu}{L+\mu}, \quad y_4=\frac{2}{L+\mu},\quad y_5=1,
 \end{equation}
and adding the five inequalities together after multiplying each one by the corresponding multiplier.

The result is the following inequality (as may be verified directly):

\begin{equation}
\begin{array}{rcl}
f_1-f_* &\leq & \left(\frac{L-\mu}{L+\mu}\right)^2 (f_0-f_*)-\frac{\mu L (L+3\mu)}{2(L+\mu)^2} \norm{\mathbf{x}_0-\frac{L+\mu}{L+3\mu}\mathbf{x}_1-\frac{2\mu}{L+3\mu} {\mathbf{x}_*-}\frac{3L+\mu}{L^2+3\mu L}\mathbf{g}_0-\frac{L+\mu}{L^2+3\mu L}\mathbf{g}_1}^2 \\
&&-\frac{2L\mu^2}{L^2+2L\mu-3\mu^2}\norm{\mathbf{x}_1-\mathbf{x}_*-\frac{(L-\mu)^2}{2\mu L(L+\mu)}\mathbf{g}_0-\frac{L+\mu}{2\mu L} \mathbf{g}_1}^2.
\end{array}
\label{eq:key inequality}
\end{equation}
Since the last two right-hand-side terms are nonpositive, we obtain:
\[
f_1-f_*\leq \left(\frac{L-\mu}{L+\mu}\right)^2 (f_0-f_*).
\]
Since $\mathbf{x}_0$ was arbitrary,
this  completes the proof of Theorem \ref{th:main}. \qed

\subsection{Remarks on the proof of Theorem \ref{th:main}.}
%A few remarks on Theorem \ref{th:main} and its proof:

\begin{itemize}
\item
First, note that we have proven a bit more than what is stated in Theorem \ref{th:main}.
Indeed, the result in  Theorem \ref{th:main} holds for any iterative method that satisfies the five inequalities used in its proof.

\item
Although the proof of Theorem \ref{th:main} is easy to verify, it
is not apparent how the multipliers $y_1,\ldots, y_5$ in \eqref{eq:multipliers} were obtained.
This was in fact done via preliminary computation{s}, and subsequently guessing the values in \eqref{eq:multipliers},
through the following steps:

\begin{enumerate}
\item
The SDP performance estimation problem
\eqref{eq:pep sdp} with $N=1$ was solved numerically for various  values of the parameters {$\mu$ , $L$ and $R$} --- actually, the values of $L$ and $R$ can safely be fixed to some positive constants using appropriate scaling arguments (see e.g.,~\cite[Section 3.5]{PEP2016smooth} for a related discussion).
\item
The optimal values
of the dual SDP multipliers of the constraints corresponding to the five inequalities in the proof gave the guesses for the
correct values $y_1,\ldots, y_5$ as stated in in \eqref{eq:multipliers}.
\item
 Finally the correctness of the guess was verified directly (by symbolic computation and by hand).
\end{enumerate}

\item
The key inequality \eqref{eq:key inequality} may be rewritten in another, more symmetric way%, under the assumption (without loss of generality) that $\mathbf{x}_* = 0$:
\[
(f_1 - f_*) \le (f_0 - f_*) \left(\frac{1-\kappa}{1+\kappa}\right)^2 - \frac{\mu}{4} \left( \frac{{\norm{\mathbf{s}_1}^2}}{1+\sqrt{\kappa}} + \frac{{\norm{\mathbf{s}_2}^2}}{1-\sqrt{\kappa}} \right),
\]
where  $\kappa = \mu/L$ is the condition number (between $0$ and $1$) and slack vectors $\mathbf{s}_1$ and $\mathbf{s}_2$ are
\begin{eqnarray*}
{\mathbf{s}_1} &=& -\frac{(1+\sqrt{\kappa})^2}{1+\kappa}  \left(\mathbf{x}_0 -\mathbf{x}_* - \mathbf{g}_0/\sqrt{L\mu} \right) + \left( \mathbf{x}_1 - \mathbf{x}_* + \mathbf{g}_1/\sqrt{L\mu} \right) \\
{\mathbf{s}_2} &=&\ \ \ \frac{(1-\sqrt{\kappa})^2}{1+\kappa}  \left( \mathbf{x}_0 - \mathbf{x}_* + \mathbf{g}_0/\sqrt{L\mu} \right) - \left(\mathbf{x}_1 - \mathbf{x}_* - \mathbf{g}_1/\sqrt{L\mu} \right). \\
\end{eqnarray*}
Note that the four expressions  $ \mathbf{x}_i  - \mathbf{x}_* \pm  \mathbf{g}_i/\sqrt{L\mu}$ expressions are invariant under dilation of $f$, and that cases of equality in \eqref{eq:key inequality} simply correspond to equalities $\mathbf{s}_1=\mathbf{s}_2=0$.

%We presented both forms of the inequality in this paper, since the explicit dependence on $\mathbf{x}_*$ in \eqref{eq:key inequality}
%may be of interest.

\item
It is interesting to note that the known proof of Theorem \ref{th:main} for the quadratic case only requires  the so-called Kantorovich inequality,
that may be stated as follows.

\begin{theorem}[Kantorovch inequality; see e.g.\ Lemma 3.1 in \cite{bertsekas1999nonlinear}]
Let $Q$ be a symmetric positive definite $n\times n$ matrix with smallest and largest eigenvalues $\mu >0$ and $L\ge \mu$ respectively.
Then, for any unit vector $\mathbf{x} \in \mathbb{R}^n$, one has:
\[
\left(\mathbf{x}^\tsp Q \mathbf{x}\right)\left(\mathbf{x}^\tsp Q^{-1} \mathbf{x}\right) \le \frac{(\mu+L)^2}{4\mu L}.
\]
\end{theorem}
Thus, the inequality \eqref{eq:key inequality} replaces the Kantorovich inequality in the proof of Theorem \ref{th:main} for non-quadratic
$f \in \mathcal{F}_{\mu,L}(\mathbb{R}^n)$.

\item \modFG{Finally, we note that this proof can be modified very easily to handle the case of the fixed-step gradient method that was mentioned in Theorem~\ref{thNesterov}. Indeed, observe that the proof aggregates the fourth and fifth inequalities with multipliers $y_4=\frac{2}{L+\mu}$ and $y_5=1$, which leads to the combined inequality \[ \frac{2}{L+\mu}
\bigl( -\mathbf{g}_{0}^{\tsp} \mathbf{g}_{1} \bigr) +
\mathbf{g}_{1}^{\tsp} (\mathbf{x}_0-\mathbf{x}_1)\geq 0
\quad \Leftrightarrow \quad \mathbf{g}_{1}^{\tsp}  (\mathbf{x}_0 - \frac{2}{L+\mu} \mathbf{g}_0 - \mathbf{x}_1 ) \ge 0 \;. \] Now note that the gradient method with fixed step $\gamma = \frac{2}{L+\mu}$ satisfies this combined inequality (since the second factor in the left-hand side becomes zero), and hence the rest of the proof establishes the same rate for this method as for the gradient descent with exact line search.
\begin{theorem}
Let $f\in \mathcal{F}_{\mu,L}(\mathbb{R}^n)$, $\mathbf{x}_*$ a global minimizer of $f$ on $\mathbb{R}^n$, and $f_* = f(\mathbf{x}_*)$.
Each iteration of the gradient method with {fixed step length} $\gamma=\frac{2}{\mu+L}$ satisfies
\begin{equation*}
{f(\mathbf{x}_{i+1})} - f_*\leq \left(\frac{L-\mu}{L+\mu}\right)^2 \left({f(\mathbf{x}_i)}-f_*\right) \quad i = 0,1,\ldots
\end{equation*}
\end{theorem}
Finally, note that Example~\ref{eg:quad example} also establishes that this rate is tight. \modFGR{Hence we have the relatively surprising fact that, when looking at the worst-case convergence rate of the objective function accuracy, performing exact line-search is not better than using a well-chosen fixed step length.}
%Interestingly, the proof of Theorem 1.2 can also be applied on the fixed-step gradient method with step size 2 , by replacing inequalities 4 and 5 on page 6 with
% L+μ
%′⊤ 2
%and taking y4′ = 1. Example 1.3 can then be used to provide the lower bound necessary
%for establishing the tightness of the resulting bound.
}\end{itemize}

\section{Extension to `noisy' gradient descent with exact line search}
Theorem \ref{th:main} may be generalized to what we will call \emph{noisy gradient descent method with exact linear search}; see e.g.\ \cite[p.59]{bertsekas1999nonlinear} where it is called gradient descent method with (relative) error.
Here the search direction at iteration $i$, say $\mathbf{d}_i$, satisfies
\begin{equation}
\label{eq:error}
\|-\nabla f(\mathbf{x}_i)-\mathbf{d}_i\| \le \varepsilon\|\nabla f(\mathbf{x}_i)  \| \quad i = 0,1,\ldots,
\end{equation}
where {$0\le \varepsilon < 1$} is some given \emph{relative} tolerance on the deviation from the negative gradient. Note that the algorithm cannot be guaranteed to converge as soon as $\varepsilon\geq 1$, since $\mathbf{d_i}=0$ then becomes feasible. %Thus the search direction may be viewed as a small perturbation of the negative gradient.
We recover the normal gradient descent algorithm when $\varepsilon = 0$.

In the case of more general values of $\varepsilon$, one can for example satisfy the relative error criterion
 by imposing a restriction of the type $|\sin \theta| \leq \varepsilon$ on the angle $\theta$ between search direction
   $\mathbf{d}_i$ and the current negative gradient $-\nabla f(\mathbf{x}_i)$.

Using a search direction $\mathbf{d}_i$ that satisfies \eqref{eq:error} corresponds, for example,
 to an implementation of the gradient descent method where each component of $-\nabla f(\mathbf{x}_i)$ is only calculated to
a fixed number of significant digits. It is also related to the so-called \emph{stochastic gradient descent} method that
is used in training neural networks;
 see e.g.\ \cite{neural}
and the references therein.

Thus we consider the following algorithm:

\begin{center}
{\small \fbox{
\parbox{0.7\textwidth}{
        \textbf{Noisy gradient descent method with exact line search}
  \begin{itemize}
  \item[] \textbf{Input:} $f\in\mathcal{F}_{\mu,L}(\mathbb{R}^n)$, $\mathbf{x}_0\in \mathbb{R}^n$, ${0\le \varepsilon < 1}$.\\[-0.2cm]
  \item[] \textbf{for} $i=0,1,\ldots$ \\ [-0.2cm]
  \item[]\hspace{1cm} Select any seach direction $\mathbf{d}_i$ that satisfies \eqref{eq:error}; \\ [-0.2cm]
  \item[]\hspace{1cm} {$\gamma=\text{argmin}_{\gamma \in \mathbb{R}}f\left(\mathbf{x}_i-\gamma \mathbf{d}_i\right)$} \\ [-0.2cm]
  \item[]
   \hspace{1cm} {$\mathbf{x}_{i+1}=\mathbf{x}_i-\gamma \mathbf{d}_i$} \\ [-0.2cm]
  \end{itemize}
       }}}
\end{center}

One may show the following generalization of Theorem \ref{th:main}.

\begin{theorem}
\label{th:noisy}
Let $f\in \mathcal{F}_{\mu,L}(\mathbb{R}^n)$, $\mathbf{x}_*$ a global minimizer of $f$ on $\mathbb{R}^n$, and $f_* = f(\mathbf{x}_*)$.
Given a relative tolerance $\varepsilon$, each iteration of the noisy gradient descent method with exact line search satisfies
\begin{equation}
\label{eq:boundnoisy}
{f(\mathbf{x}_{i+1})} - f_*\leq \modAT{\left(\frac{1-\kappa_{\varepsilon}}{1+\kappa_{\varepsilon}}\right)^2} ({f(\mathbf{x}_i)}-f_*) \quad i = 0,1,\ldots
\end{equation}
where \modAT{$\kappa_{\varepsilon}= \frac{\mu}{L} \frac{(1-\varepsilon)}{(1+\varepsilon)}$.}
\end{theorem}
\modFGR{When $\varepsilon = 0$, the rate becomes $\frac{1-\kappa}{1+\kappa} = \frac{L-\mu}{L+\mu}$, which matches} exactly Theorem \ref{th:main}, and
the proof of Theorem \ref{th:noisy} is a straightforward generalization of the proof of Theorem {\ref{th:main}. The} key is again to consider a wider class of iterative methods that satisfies certain inequalities.
Here we use the inequalities:

\begin{equation}
\left.
\begin{array}{lcl}
1: &&f_{0}\geq f_1 + \inner{\mathbf{g}_1}{(\mathbf{x}_{0}-\mathbf{x}_1)} +\frac{1}{2(1-\mu/L)}\left( \frac{1}{L}\norm{\mathbf{g}_0-\mathbf{g}_1}^2+ \mu \norm{\mathbf{x}_0-\mathbf{x}_1}^2 - 2\frac{\mu}{L} \inner{(\mathbf{g}_1-\mathbf{g}_0)}{(\mathbf{x}_1-\mathbf{x}_0)}\right)\\
2: &&f_*\geq f_0 + \inner{\mathbf{g}_0}{(\mathbf{x}_*-\mathbf{x}_0)} +\frac{1}{2(1-\mu/L)}\left( \frac{1}{L}\norm{\mathbf{g}_*-\mathbf{g}_0}^2+ \mu \norm{\mathbf{x}_*-\mathbf{x}_0}^2 - 2\frac{\mu}{L} \inner{(\mathbf{g}_0-\mathbf{g}_*)}{(\mathbf{x}_0-\mathbf{x}_*)}\right)\\
3: &&f_*\geq f_1 + \inner{\mathbf{g}_1}{(\mathbf{x}_*-\mathbf{x}_1)} +\frac{1}{2(1-\mu/L)}\left( \frac{1}{L}\norm{\mathbf{g}_*-\mathbf{g}_1}^2+ \mu \norm{\mathbf{x}_*-\mathbf{x}_1}^2 - 2\frac{\mu}{L} \inner{(\mathbf{g}_1-\mathbf{g}_*)}{(\mathbf{x}_1-\mathbf{x}_*)}\right)\\
4: && 0 \geq \mathbf{g}_{1}^\tsp (\mathbf{x}_{1}-\mathbf{x}_0)\\
5: && 0 \ge \mathbf{g}_0^\tsp \mathbf{g}_1- \varepsilon \|\mathbf{g}_0\|\|\mathbf{g}_1\|.\\
\end{array}
\right\}
\label{noisy constraints}
\end{equation}
The first four inequalities are the same as before, and the fifth is satisfied by the iterates of the noisy gradient descent with exact line search. Indeed, in the first iteration one has:
\begin{eqnarray*}
0 &=& \mathbf{d}_0^T \frac{\mathbf{g}_1}{\|\mathbf{g}_1\|} \quad\quad\quad\quad\quad\quad\quad\quad\quad\ \mbox{ (exact line search)}\\
  &=& (\mathbf{d}_0 + \mathbf{g}_0)^\tsp \frac{\mathbf{g}_1}{\|\mathbf{g}_1\|} - \frac{\mathbf{g}_0^\tsp \mathbf{g}_1}{\|\mathbf{g}_1\|} \\
  &\le& \varepsilon \|\mathbf{g}_0\| - \frac{\mathbf{g}_0^\tsp \mathbf{g}_1}{\|\mathbf{g}_1\|} \quad\quad\quad\quad\quad\quad\quad \mbox{(by Cauchy-Schwartz and \eqref{eq:error}).} \\
\end{eqnarray*}

We rewrite the fifth inequality as the equivalent linear matrix inequality:
\begin{equation}
\label{noisy lmi}
\left(
     \begin{array}{cc}
      \varepsilon \|\mathbf{g}_0\|^2 & \mathbf{g}_0^\tsp \mathbf{g}_1\\
\mathbf{g}_0^\tsp \mathbf{g}_1 & \varepsilon \|\mathbf{g}_1\|^2 \\
     \end{array}
   \right)
\succeq 0.
\end{equation}

We first aggregate the first four inequalities in \eqref{noisy constraints} by adding them together after multiplication by the respective multipliers:\modAT{
\[
 y_1= \rho_\varepsilon,\quad y_2=2\kappa_\varepsilon \frac{1-\kappa_\varepsilon}{(1+\kappa_\varepsilon)^2}, \quad y_3=\frac{2\kappa_\varepsilon}{1+\kappa_\varepsilon}, \quad y_4=1,
\]}
where
{ $L_\varepsilon=(1+\varepsilon)L$, $\mu_\varepsilon=(1-\varepsilon)\mu$, $\kappa_\varepsilon=\frac{\mu_\varepsilon}{L_\varepsilon}$ and $\rho_\varepsilon=\frac{1-\kappa_\varepsilon}{1+\kappa_\varepsilon}$}.

Next we define a positive semidefinite matrix multiplier for the linear matrix inequality \eqref{noisy lmi}, namely
 \begin{equation}
\label{noisy lmi multiplier}\modAT{
 \begin{pmatrix}
a\rho_\varepsilon & -a \\
-a & \frac{a}{\rho_\varepsilon}
\end{pmatrix}\succeq 0,}
\end{equation}
\modAT{with $a=\frac{1}{L_\varepsilon+\mu_\varepsilon}$},
 and add \modFG{nonnegativity of} the inner product
 between the left-hand-side of \eqref{noisy lmi} and the multiplier matrix \eqref{noisy lmi multiplier} to the aggregated constraints. It can now be checked that the resulting expression is the following (slight) generalization of \eqref{eq:key inequality} \modAT{
 \begin{align*}
f_1-f_*\leq \rho_\varepsilon^2 (f_0-f_*)&-\frac{L\mu (L_\varepsilon-\mu_\varepsilon)(L_\varepsilon+3\mu_\varepsilon)}{2(L-\mu)(L_\varepsilon+\mu_\varepsilon)^2} \norm{\mathbf{x}_0+\alpha_1 \mathbf{x}_1-(1+\alpha_1) \mathbf{x}_*+\alpha_2 \mathbf{g}_0+\alpha_3 \mathbf{g}_1}^2\\
&-\frac{2L\mu\mu_\varepsilon}{(L-\mu)(L_\varepsilon+3\mu_\varepsilon)}\norm{\mathbf{x}_1-\mathbf{x}_*+\alpha_4 \mathbf{g}_0+\alpha_5 \mathbf{g}_1}^2,
\end{align*}
with the appropriate coefficients
\[
\alpha_1=-\frac{L_\varepsilon+\mu_\varepsilon}{L_\varepsilon+3\mu_\varepsilon}, \; \alpha_2=-\frac{4L-L_\varepsilon+\mu_\varepsilon}{L(L_\varepsilon+3\mu_\varepsilon)}, \; \alpha_3=\frac{(L_\varepsilon+\mu_\varepsilon)(-4L+3L_\varepsilon+\mu_\varepsilon)}{L(L_\varepsilon-\mu_\varepsilon)(L_\varepsilon+3\mu_\varepsilon)}, \; \alpha_4=-\frac{(L-\mu)(L_\varepsilon-\mu_\varepsilon)}{2L\mu(L_\varepsilon+\mu_\varepsilon)},
\]
and $\alpha_5=-\frac{L+\mu}{2L\mu}$.
 This completes the proof. \qed}

\par

\modFG{To conclude this section, the following example, based on the same quadratic function as Example\ref{eg:quad example}, shows that our bound \eqref{eq:boundnoisy} for the noisy gradient descent is also tight.

\begin{example}\label{eg:quad example noisy}
Consider the same quadratic function as in Example~\ref{eg:quad example}:
\[
f(\mathbf{x}) = \frac{1}{2}\sum_{i=1}^n \lambda_i x_i^2
\quad \text{ where } \quad
0 < \mu = \lambda_1 \le \lambda_2 \le \ldots \le \lambda_{n} = L.
\]
Let $\theta$ be an angle satisfying $0 \le \theta < \frac{\pi}{2}$. Consider the noisy gradient descent method where direction $\mathbf{d}_0$ is obtained by performing a clockwise 2D-rotation with angle $\theta$ on the first and last coordinates of the gradient $\nabla f(\mathbf{x}_0)$. As mentioned above, this satisfies our definition with relative tolerance $\varepsilon = \sin \theta$. Define now the starting point
\[
\mathbf{x}_0 = \left(\frac1\mu   , 0, \ldots, 0, \frac1L \modFGR{\sqrt{\frac{1-\varepsilon}{1+\varepsilon}} } \right)^\tsp .
\]
Tedious but straightforward computations  show that
\[ \mathbf{x}_1 = \left(\frac{1-\kappa_{\varepsilon}}{1+\kappa_{\varepsilon}}\right)   \left(\frac1\mu  , 0, \ldots, 0, -\frac1L \modFGR{\sqrt{\frac{1-\varepsilon}{1+\varepsilon}} } \right)^\tsp \quad \text{ where } \kappa_{\varepsilon}= \frac{\mu}{L} \frac{(1-\varepsilon)}{(1+\varepsilon)}. \]
Moreover, if one chooses $\mathbf{d}_1$ by rotating the second gradient $\nabla f(\mathbf{x}_1)$ by the same angle $\theta$ in the counterclockwise direction, one obtains
\[ \mathbf{x}_2 = \left(\frac{1-\kappa_{\varepsilon}}{1+\kappa_{\varepsilon}}\right)^2 \left(\frac1\mu   , 0, \ldots, 0, \frac1L \modFGR{\sqrt{\frac{1-\varepsilon}{1+\varepsilon}} } \right)^\tsp = \left(\frac{1-\kappa_{\varepsilon}}{1+\kappa_{\varepsilon}}\right)^2 \mathbf{x}_0. \]
A similar reasoning for the next iterates, alternating clockwise and counterclockwise rotations, shows that
\[
\mathbf{x}_{2i} = \left(\frac{1-\kappa_{\varepsilon}}{1+\kappa_{\varepsilon}}\right)^{2i}x_0, \;\;\; \mathbf{x}_{2i+1} = \left(\frac{1-\kappa_{\varepsilon}}{1+\kappa_{\varepsilon}}\right)^{2i}x_1 \; \text{ for all } i=0,1,\ldots
\]
and hence we have that equality
\[
{f(\mathbf{x}_{i+1})} - f_* = \left(\frac{1-\kappa_{\varepsilon}}{1+\kappa_{\varepsilon}}\right)^2 ({f(\mathbf{x}_i)}-f_*) \quad i = 0,1,\ldots
\]
holds as announced. Figure~\ref{fig:quad_ex_noisy} displays a few iterates, and can be compared to Figure~\ref{fig:quad_ex}.\qed
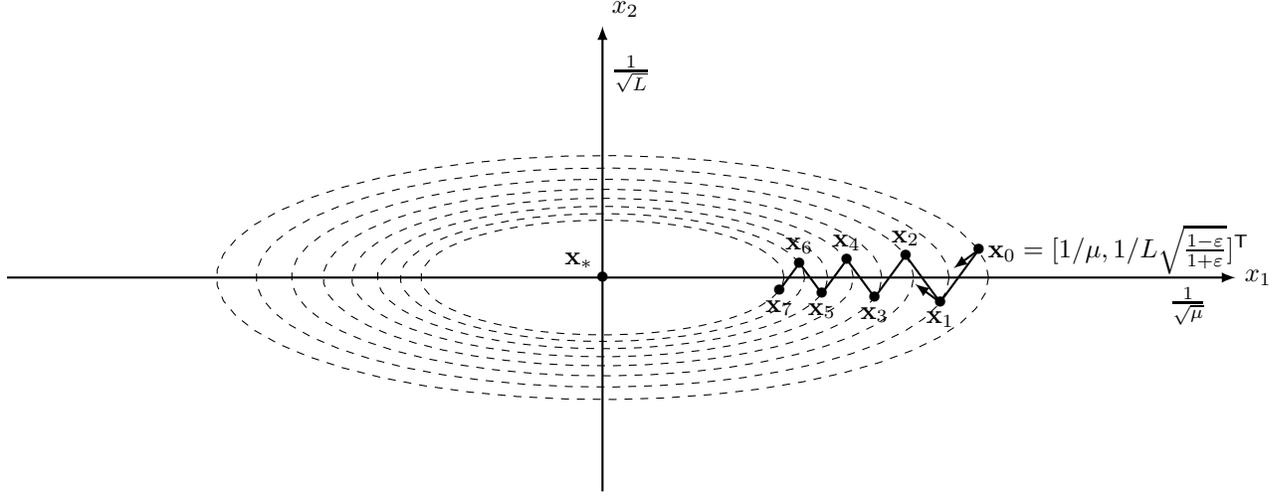
\begin{figure}
\begin{center}
\begin{tikzpicture}[yscale=1,xscale=1]
\tikzset{fleche/.style={->, >=latex, thick},
	    flecheB/.style={<->, >=latex, very thick},
 	    droite/.style={thick, dashed}};
 	
\pgfmathsetmacro{\xx}{5};
\pgfmathsetmacro{\yy}{0.3669};
\pgfmathsetmacro{\rhoExB}{0.89781};
\pgfmathsetmacro{\fxs}{(1/1+1/.1)/2};
\pgfmathsetmacro{\lAxis}{sqrt(2*\fxs/.2)}
\pgfmathsetmacro{\sAxis}{sqrt(2*\fxs/2)};

\pgfmathsetmacro{\nbEllipses}{8};%Number of ellipses to show
\pgfmathsetmacro{\scaleEllipse}{sqrt(\xx^2/10+\yy^2)};

\foreach \i in {1,...,\nbEllipses}
{
\draw [dashed, color=black, domain=-0:360, samples=50] plot({\scaleEllipse*cos(\x)*sqrt(10)*\rhoExB^(\i-1)}, {\scaleEllipse*sin(\x)*sqrt(1)*\rhoExB^(\i-1)});
}
\node (XS) at (0,0) {$\bullet$};
\draw (XS) node[above left] {$\mathbf{x}_*$};
\node (X0) at (\xx,\yy) {$\bullet$};
\draw (X0) node[right ] {${\mathbf{x}_0} = [1/\mu,1/L\sqrt{\frac{1-\varepsilon}{1+\varepsilon}}]^\tsp$};
\node (X1) at (\rhoExB*\xx,-\rhoExB*\yy) {$\bullet$};
\draw (X1) node[below ] {$\mathbf{x}_1$};
\node (X2) at (\rhoExB^2*\xx,\rhoExB^2*\yy) {$\bullet$};
\draw (X2) node[above ] {$\mathbf{x}_2$};
\node (X3) at (\rhoExB^3*\xx,-\rhoExB^3*\yy) {$\bullet$};
\draw (X3) node[below ] {$\mathbf{x}_3$};
\node (X4) at (\rhoExB^4*\xx,\rhoExB^4*\yy) {$\bullet$};
\draw (X4) node[above ] {$\mathbf{x}_4$};
\node (X5) at (\rhoExB^5*\xx,-\rhoExB^5*\yy) {$\bullet$};
\draw (X5) node[below ] {$\mathbf{x}_5$};
\node (X6) at (\rhoExB^6*\xx,\rhoExB^6*\yy) {$\bullet$};
\draw (X6) node[above ] {$\mathbf{x}_6$};
\node (X7) at (\rhoExB^7*\xx,-\rhoExB^7*\yy) {$\bullet$};
\draw (X7) node[below ] {$\mathbf{x}_7$};
\draw[thick] (X0.center)--(X1.center)--(X2.center)--(X3.center)--(X4.center)--(X5.center) -- (X6.center) -- (X7.center);
\draw[fleche] (X0.center)--(\xx-.2*\xx/3,\yy-2*\yy/3);
\draw[fleche] (X1.center)--(\rhoExB*\xx-.2*\xx/3,-\rhoExB*\yy+2*\yy/3);
\draw [fleche] (-\lAxis-.5,0) -- (\lAxis+1,0);
\draw [fleche] (0,-\sAxis-.5) -- (0,\sAxis+1);
\draw (\lAxis+1,0) node[ right] {$x_1$};
\draw (0,\sAxis+1) node[above right] {$x_2$};
\draw (0,\sAxis) node[above right] {$\frac{1}{\sqrt{L}}$};
\draw (\lAxis,0) node[below right] {$\frac{1}{\sqrt{\mu}}$};
\end{tikzpicture}
\caption{\modFGR{Illustration Example \ref{eg:quad example noisy}  for $n=2$ and $\varepsilon=0.3$ (small arrows indicate direction of negative gradient). \label{fig:quad_ex_noisy}}}
\end{center}
\end{figure}

\end{example}
}

\section{Concluding remarks}
 The main results of this paper are the exact convergence rates of the gradient descent method with exact line search \modFG{and its noisy variant} %(and a bound on the worst-case convergence of its noisy variant)
for strongly convex functions with Lipschitz continuous gradients. The computer-assisted technique of proof is also of independent interest, and demonstrates the importance of the SDP performance estimation problems (PEPs) introduced in \cite{drori2014}.

Indeed, to obtain our proof of Theorem \ref{th:noisy}, the following SDP PEP was solved numerically for various fixed values of $R$, $\mu$ and $L$:
\[
\max f_1 - f_* \; \mbox{ subject to } \eqref{noisy constraints} \mbox{ and } f_0 - f_* \le R.
\]
It was observed that, for each set of values, the optimal value of the SDP corresponded exactly to the bound in Theorem~\ref{th:noisy} (actually, for homogeneity reasons, $L$ and $R$ could be fixed and only $\mu$ needed to vary). Based on this, a rigorous proof Theorem~\ref{th:noisy} could be given by guessing the correct values of the dual SDP multipliers as functions of $\mu$, $L$ and $R$, and then verifying the guess through an explicit computation.

We believe this type of computer-assisted proof could prove useful in the analysis of more methods where exact line search is used (see for example ~\cite{PEP2016composite} which studies conditional gradient methods).

PEPs have been used by now to study worst-case convergence rates of several first-order optimization methods \cite{drori2014,PEP2016smooth,PEP2016composite}. This paper differs in an important aspect:  the performance estimation problem considered actually characterizes a whole class of methods that contains the method of interest (gradient descent with exact line search) as well as many other methods. This relaxation in principle only provides an upper bound on the worst-case of gradient descent, and it is the fact that Example~\ref{eg:quad example} matches this bound that allows us to conclude with a tight result.

The reason we could not solve the peformance estimation problem for the gradient descent method itself is that equation \eqref{eq:negative gradient}, which essentially states that the step $\mathbf{x}_{i+1}-\mathbf{x}_i$ is parallel to the gradient $\nabla f(\mathbf{x}_i)$, cannot be formulated as \modATR{a convex} constraint \modATR{in the SDP formulation}. The main obstruction appears to be that requiring that two vectors are parallel is a nonconvex constraint, even when working with their inner products\footnote{One such nonconvex formulation would be $\mathbf{g}_{i}^{\tsp} (\mathbf{x}_i-\mathbf{x}_{i+1}) = \|\mathbf{g}_{i}\| \|\mathbf{x}_i-\mathbf{x}_{i+1} \|$.}. Instead, our convex formulation enforces  that those two vectors are both orthogonal to a third one, the next gradient $\nabla f(\mathbf{x}_{i+1})$.

% The difference with the current paper, is that it was often only possible to obtain upper bounds on the convergence rate.
%In this paper the PEPs with exact line search gave tight upper bounds, and based on the dual SDP solutions it was possible to construct
%rigorous proofs of the convergence rate.
%

%as for example , which can also be modelled using semidefinite programming~\cite{PEP2016composite})}.

\subsubsection*{Acknowledgements}
The authors would like to thank Simon Lacoste-Julien for bringing Theorem 2.1.15 in \cite{nesterov2004intro} to their attention, and an anonymous referee
\modFG{for valuable suggestions that include the last remark in Section 4.1}.


\begin{thebibliography}{10}

\bibitem{bertsekas1999nonlinear}
D.P. Bertsekas.
\newblock {\em Nonlinear programming}.
\newblock Athena scientific, 1999.

\bibitem{drori2014}
Y.~Drori and M.~Teboulle.
\newblock Performance of first-order methods for smooth convex minimization: a
  novel approach.
\newblock {\em Mathematical Programming}, 145(1-2):451--482, 2014.

\bibitem{luenberger2008linear}
D.G. Luenberger and Y.~Ye.
\newblock {\em Linear and nonlinear programming}.
\newblock Springer, 2008.

\bibitem{neural}
A.~Neelakantan, L.~Vilnis, Q.V. Le, I.~Sutskever, L.~Kaiser, K.~Kurach, and
  J.~Martens.
\newblock Adding gradient noise improves learning for very deep networks.
\newblock {\em arXiv}, 1511.06807v1, 2015.

\bibitem{nemirovskinotes1999}
A.~Nemirovski.
\newblock Optimization {I}{I}: Numerical methods for nonlinear continuous
  optimization.
\newblock Lecture notes, 1999.
\newblock Available from:
  \url{http://www2.isye.gatech.edu/~nemirovs/Lect_OptII.pdf}.

\bibitem{nemirovsky1983problem}
A.~Nemirovski and D.B. Yudin.
\newblock {\em Problem complexity and method efficiency in optimization}.
\newblock Wiley, 1983.

\bibitem{nesterov2004intro}
Yu. Nesterov.
\newblock {\em Introductory lectures on convex optimization : a basic course}.
\newblock Applied optimization. Kluwer Academic Publ., 2004.

\bibitem{nocedal2006numerical}
J.~Nocedal and S.~Wright.
\newblock {\em Numerical optimization}.
\newblock Springer Science \& Business Media, 2006.

\bibitem{polyak1987intro}
B.T. Polyak.
\newblock {\em Introduction to optimization}.
\newblock Optimization Software, New York, 1987.

\bibitem{PEP2016composite}
A.B. Taylor, J.M. Hendrickx, and F.~Glineur.
\newblock Exact worst-case performance of first-order methods for composite
  convex optimization.
\newblock {\em arXiv}, 1512.07516, 2015.

\bibitem{PEP2016smooth}
A.B. Taylor, J.M. Hendrickx, and F.~Glineur.
\newblock Smooth strongly convex interpolation and exact worst-case performance
  of first-order methods.
\newblock {\em Mathematical Programming}, 2016.
\newblock Accepted for publication.

\end{thebibliography}
 \end{document}